\documentclass[11pt]{amsart}
\usepackage{amsmath,amsfonts,amsthm,amsopn,amssymb,mathtools,stmaryrd}
\usepackage{cite,marginnote}
\usepackage{pdfsync}

\usepackage{color,enumitem,graphicx}
\usepackage[colorlinks=true,urlcolor=blue,
citecolor=red,linkcolor=blue,linktocpage,pdfpagelabels,
bookmarksnumbered,bookmarksopen]{hyperref}
\usepackage[english]{babel}

\usepackage[left=2.75cm,right=2.75cm,top=3.4cm,bottom=3.4cm]{geometry}
\usepackage[]{xcolor}

\numberwithin{equation}{section}

\makeindex
\newtheorem{theorem}{Theorem}[section]
\newtheorem{proposition}[theorem]{Proposition}
\newtheorem{lemma}[theorem]{Lemma}
\newtheorem{remark}[theorem]{Remark}
\newtheorem{example}[theorem]{Example}
\newtheorem{corollary}[theorem]{Corollary}
\newtheorem{definition}[theorem]{Definition}

\newcommand{\bt}{\begin{theorem}}
\newcommand{\et}{\end{theorem}}
\newcommand{\bl}{\begin{lemma}}
\newcommand{\el}{\end{lemma}}
\newcommand{\bd}{\begin{definition}}
\newcommand{\ed}{\end{definition}}
\newcommand{\bc}{\begin{corollary}}
\newcommand{\ec}{\end{corollary}}
\newcommand{\bp}{\begin{proof}}
\newcommand{\ep}{\end{proof}}
\newcommand{\bx}{\begin{example}}
\newcommand{\ex}{\end{example}}
\newcommand{\bi}{\begin{exercise}}
\newcommand{\ei}{\end{exercise}}
\newcommand{\bo}{\begin{proposition}}
\newcommand{\eo}{\end{proposition}}
\newcommand{\br}{\begin{remark}}
\newcommand{\er}{\end{remark}}
\newcommand{\be}{\begin{equation}}
\newcommand{\ee}{\end{equation}}
\newcommand{\ba}{\begin{align}}
\newcommand{\ea}{\end{align}}
\newcommand{\bn}{\begin{enumerate}}
\newcommand{\en}{\end{enumerate}}
\newcommand{\bg}{\begin{align*}}
\newcommand{\eg}{\end{align*}}
\newcommand{\bcs}{\begin{cases}}
\newcommand{\ecs}{\end{cases}}
\newcommand{\bean}{\begin{eqnarray*}}
\newcommand{\eean}{\end{eqnarray*}}

\newcommand{\R}{\mathbb R}
\newcommand{\al}{\alpha}

\newcommand{\C}{\mathbb C}
\newcommand{\vs}{\vskip .25in}

\newcommand{\N}{{\mathcal N}}
\newcommand{\M}{{\mathcal M}}

\title[Ground states of coupled Schr\"{o}dinger systems]{Ground states of Schr\"{o}dinger systems \\ with Chern-Simons gauge fields}

\author[Y. H.\ Jiang]{Yahui Jiang}
\author[T. Y.\ Chen]{Taiyong Chen}
\author[J. J.\ Zhang]{Jianjun Zhang}
\author[M. Squassina]{Marco Squassina}
\author[N. Almousa]{Nouf Almousa}

\address[Y. H.\ Jiang]{\newline\indent School of Mathematics
\newline\indent
China University of Mining and Technology
\newline\indent
Xuzhou, 221116, China}
\email{\href{mailto:18843111149@163.com}{18843111149@163.com}}

\address[T. Y.\ Chen]{\newline\indent School of Mathematics
\newline\indent
China University of Mining and Technology
\newline\indent
Xuzhou, 221116, China}
\email{\href{mailto:taiyongchencumt@163.com}{taiyongchencumt@163.com}}

\address[J. J.\ Zhang]{\newline\indent College of Mathematica and Statistics
\newline\indent
Chongqing Jiaotong University
\newline\indent
Chongqing 400074, China}
\email{\href{mailto:zhangjianjun09@tsinghua.org.cn}{zhangjianjun09@tsinghua.org.cn}}

\address[M.\ Squassina]{\newline\indent College of Science
	\newline\indent
	Princess Nourah Bint Abdul Rahman University
	\newline\indent
	Saudi Arabia, Riyadh, PO Box 84428}
\email{\href{mailto:marsquassina@pnu.edu.sa}{marsquassina@pnu.edu.sa}}

\address[N. Almousa]{\newline\indent College of Science
	\newline\indent
	Princess Nourah Bint Abdul Rahman University
	\newline\indent
	Saudi Arabia, Riyadh, PO Box 84428}
\email{\href{mailto:nmalmousa@pnu.edu.sa}{nmalmousa@pnu.edu.sa}}

\address[M.\ Squassina]{\newline\indent Dipartimento di Matematica e Fisica
	\newline\indent
	Universit\`a Cattolica del Sacro Cuore
	\newline\indent
	Via dei Musei 41, Brescia, Italy}
\email{\href{mailto:marco.squassina@unicatt.it}{marco.squassina@unicatt.it}}

\thanks{(1) Corresponding author: marsquassina@pnu.edu.sa}

\thanks{(2) J. J.\ Zhang was partially supported by NSFC(No.11871123).}


\date{\today}


\begin{document}

\begin{abstract}
We are concerned with the following coupled nonlinear Schr\"{o}dinger system
\begin{align*}
\bcs
-\Delta u+u+\big(\int_{|x|}^{\infty}\frac{h(s)}{s}u^{2}(s)ds+\frac{h^{2}(|x|)}{|x|^{2}}\big)u=|u|^{2p-2}u+b|v|^{p}|u|^{p-2}u,\ \ x\in\R^{2},\\
-\Delta v+\omega v+\big(\int_{|x|}^{\infty}\frac{g(s)}{s}v^{2}(s)ds+\frac{g^{2}(|x|)}{|x|^{2}}\big)v=|v|^{2p-2}v+b|u|^{p}|v|^{p-2}v,\ \ x\in\R^{2},
\ecs
\end{align*}
where $\omega,b>0$, $p>1$. By virtue of the variational approach, we show the existence of nontrivial ground state solutions. Precisely, the system above admits a positive ground state solution if $p>3$ and $b>0$ large enough or if $p\in(2,3]$ and $b>0$ small.
\end{abstract}

\maketitle

\section{Introduction}
In this paper, we consider the following coupled Schr\"odinger equations with Chern-Simons gauge fields
\begin{align}
\label{question1}
\bcs
-\Delta u+u+\big(\int_{|x|}^{\infty}\frac{h(s)}{s}u^{2}(s)ds+\frac{h^{2}(|x|)}{|x|^{2}}\big)u=|u|^{2p-2}u+b|v|^{p}|u|^{p-2}u,\ \ x\in\R^{2},\\
-\Delta v+\omega v+\big(\int_{|x|}^{\infty}\frac{g(s)}{s}v^{2}(s)ds+\frac{g^{2}(|x|)}{|x|^{2}}\big)v=|v|^{2p-2}v+b|u|^{p}|v|^{p-2}v,\ \ x\in\R^{2},
\ecs
\end{align}
where $\omega>0,b>0$, $p>1$ and
$$
h(s)=\int_{0}^{s}\frac{r}{2}u^{2}(r)dr,\qquad  g(s)=\int_{0}^{s}\frac{r}{2}v^{2}(r)dr.
$$

When $b=0$, then system (\ref{question1}) is uncoupled and it reduces to two equations of the same type. In recently years, a single nonlinear Schr\"odinger equations coupled with the Chern-Simons gauge field as follows has received much attention
\be
\label{question2}
\bcs
iD_{0}\phi+(D_{1}D_{1}+D_{2}D_{2})\phi=-f(\phi),\\
\partial_{0}A_{1}-\partial_{1}A_{0}=-{\rm Im}(\bar{\phi}D_{2}\phi),\\
\partial_{0}A_{2}-\partial_{2}A_{0}={\rm Im}(\bar{\phi}D_{1}\phi),\\
\partial_{1}A_{2}-\partial_{2}A_{1}=-\frac{1}{2}|\phi|^{2},
\ecs
\ee
where $i$ denotes the imaginary unit, $\partial_{0}=\frac{\partial}{\partial t}$, $\partial_{1}=\frac{\partial}{\partial x_{1}}$, $\partial_{2}=\frac{\partial}{\partial x_{2}}$, $(t,x_{1},x_{2})\in\mathbb{R}^{1+2}$, $\phi:\mathbb{R}^{1+2}\rightarrow\C$ is the complex scalar field, $A_{\mu}:\mathbb{R}^{1+2}\rightarrow\mathbb{C}$ is the gauge field and $D_{\mu}=-\partial_{\mu}+iA_{\mu}$ is the covariant derivative for $\mu=0,1,2$. The Chern-Simons-Schr\"odinger system consists of Schr\"odinger equations augmented by the gauge field, which was first proposed and studied in \cite{R,S}. The model was proposed to study vortex solutions, which carry both electric and magnetic charges. This feature of the model is important for the study of the high-temperature superconductor, fractional quantum Hall effect and Aharovnov-Bohm scattering. For more details about system (\ref{question2}), we refer the readers to \cite{G,CR,C}. System (\ref{question2}) is invariant under gauge transformation
\begin{align*}
\phi\rightarrow\phi e^{i\chi},\ \ A_{\mu}=A_{\mu}-\partial_{\mu}\chi,
\end{align*}
for any arbitrary $C^{\infty}$ function $\chi$. \

Byeon, Hun and Seok \cite{J-H-J} investigated the existence of standing wave solutions for system \eqref{question2} with power type nonlinearity, that is $f(u)=\lambda|u|^{p-2}u$ with $p>2$ and $\lambda>0$. By using the ansatz
\begin{align*}
\bcs
\phi(t,x)=u(|x|)e^{i\omega t},\quad A_{0}(t,x)=k(|x|),\\
A_{1}(t,x)=\frac{x_{2}}{|x|^{2}}h(|x|),\quad A_{2}(t,x)=-\frac{x_{1
}}{|x|^{2}}h(|x|).
\ecs
\end{align*}
Byeon et al. got the following nonlocal semilinear
elliptic equation
\begin{align}
\label{eq1}
-\Delta u+(\omega+\xi)u+\Big(\int_{|x|}^{\infty}\frac{h(s)}{s}u^{2}(s)ds+\frac{h^{2}(|x|)}{|x|^{2}}\Big)u=f(u),\ \ \ \ x\in\mathbb{R}^{2},
\end{align}
where $\xi$ is a constant and $h(s)$ is defined as above.
Byeon et al. showed that the existence and nonexistence of positive solutions for (\ref{eq1}) were established depending on the range of $p>2$ and $\lambda>0$. For the special case $p=4$, there exist solutions if $\lambda>1$. It seems hard to obtain the boundedness of Palais-Smale sequence when $p\in(4,6)$. They constructed a Nehari-Pohozaev manifold to obtain the boundedness of Palais-Smale sequence. For $p\in(2,4)$, Pomponio and Ruiz \cite{P-R} proved the existence and nonexistence of positive solutions for (\ref{eq1}) under the different range of the $\omega$.  A series of existence and nonexistence results of solutions for (\ref{eq1}) has been researched in \cite{JHJ,P-R,H-H,AD,C-F,Jian}.\
There has been increasing interest in studying the existence, multiplicity and the concentration behavior of the solutions for the problems above mentioned.

Problem (\ref{question1}) is a nonlocal problem due to the appearance of the term
$
\int_{|x|}^{\infty}\frac{h(s)}{s}u^{2}(s)ds,
$
which indicates that (\ref{question1}) is not a pointwise identity. This causes some mathematical difficulties which make the study of such a problem particularly interesting. System (\ref{question1}) is quite different from the following local scalar field system
\be
\label{question3}
\bcs
-\Delta u+u=|u|^{2q-2}+b|v|^{q}|u|^{q-2}u,\ \ \mbox{in}\ \mathbb{R}^{N},\\
-\Delta v+\omega^{2}v=|v|^{2q-2}+b|u|^{q}|v|^{q-2}v,\ \ \mbox{in}\ \mathbb{R}^{N},
\ecs
\ee
for $\omega>0$, $b\in\mathbb{R}$, $q\in(2,2^{*})$, which does not depend on the nonlocal term any more. The coupled nonlinear Schr\"odinger system (\ref{question3}) has attracted considerable attention in the past fifteen years. Maia et al.\cite{M-M-P} by using the variational methods and the ideas of Rabinowitz \cite{Rab} investigated the existence of positive ground state solutions for system (\ref{question3}) under sufficient conditions on parameter $b$ and $\omega$. For more progress in this aspect, we refer to \cite{A-E,AE,B-W,Lin,Sirakov} and the references therein.
\

The energy functional for system (\ref{question1}) $I:E\rightarrow\R$ is defined by
\begin{align*}
I(u,v)=\frac{1}{2}\|(u,v)\|_{E}^{2}+\frac{1}{2}\Big(B(u)+B(v)\Big)-\frac{1}{2p}F(u,v),
\end{align*}
where
\begin{align*}
&B(u)=\int_{\R^{2}}\frac{u^{2}}{|x|^{2}}\Big(\int_{0}^{|x|}\frac{s}{2}u^{2}(s)ds\Big)^{2}dx,\\
&F(u,v)=\int_{\R^{2}}\big(u^{2p}+v^{2p}+2b|uv|^{p}\big)dx.
\end{align*}
Here $E$ denotes the subspace of radially symmetric functions in $H_{r}^{1}(\R^{2})\times H_{r}^{1}(\R^{2})$ with the norm
\begin{align*}
\|(u,v)\|_{E}^{2}:=\int_{\R^{2}}\big(|\nabla u|^{2}+u^{2}+|\nabla v|^{2}+\omega v^{2}\big)dx.
\end{align*}
The functional $I$ is of class $C^{1}(E)$, its critical point $(u,v)$ is a weak solution of (\ref{question1}) and by standard regularity theory is a classical solution.\

Motivated by \cite{M-M-P}, we try to study the existence of positive ground state solutions for coupled Schr\"odinger equations with Chern-Simons gauge fields (\ref{question1}) with suitable conditions on $\omega,b$.

One of the main difficulties is the boundedness of Palais-Smale sequences if we try to use directly by the mountain pass theorem to get the critical points of $I$ in $E$. For $p\geq3$, it is standard to show that Palais-Smale condition holds for $I$. For $p\in(2,3)$, the functional $I$ has the mountain-pass geometry. However, it seems hard to prove the Palais-Smale condition holds for the functional $I$. Motivated by \cite{J-H-J}, by using a constrained minimizer on Nehari-Pohozaev manifold we circumvent this obstacle.\

Another problem is the existence of positive ground state for system (\ref{question1}), i.e., a minimal action solution $(u,v)$ with both $u>0$, $v>0$ nontrivial. We point out that the system (\ref{question1}) also possesses a trivial solution $(0,0)$ and semi-trivial solutions of type $(u,0)$ or $(0,v)$. A solution $(u,v)$ of (\ref{question1}) is nontrivial if $u\not\equiv0$ and $v\not\equiv0$. Here we overcome this obstacle by energy estimation.

We now state the main result of the paper. The constants $b_{1},b_{2},b_{3},b_{4},b_{\delta}$ involved in the statement
depends on the ground state of the single equation. We will give the corresponding expressions in Section 3.

\bt
\label{Theorem 1.1.}{\it Assume that one of the following conditions holds\\
$(i)$ $p\in(2,3]$ and $b\in(0,b_{\delta})$ sufficiently small,\\
$(ii)$ $p\in(3,3+\sqrt{6})$ and $b>\max\{b_{1},b_{2}\}$,\\
$(iii)$ $p\in[3+\sqrt{6},\infty)$ and $b>\max\{b_{3},b_{4}\},$\\
then the system (\ref{question1}) admits a positive vector ground state.}
\et

Additionally, we prove also the following nonexistence result.

\bt\label{Theorem 1.2.}
{There exist $\widetilde{b}>0$ sufficiently small and $\widetilde{\omega}>0$ sufficiently large such that $b\in(0,\widetilde{b})$ and $\omega>\widetilde{\omega}$, then the system \eqref{question1} has only trivial solution if $p\in(1,2]$.}
\et

Compared with the case of $p>2$, it seems that the case of $p\in(1,2]$ becomes more complicated and we will consider it in a forthcoming paper.
The rest of this paper is organized as follows. In Section 2, we present some notations and preliminary results and prove the nonexistence result Theorem \ref{Theorem 1.2.}. Then we give the proof of the existence of a positive ground state in Theorem \ref{Theorem 1.1.}. The item $(ii)$ and $(iii)$ are proved in Section 3, and the item $(i)$ is proved in Section 4.

\section{Preliminaries}
\noindent To prove the main results, we use the following notations:\\
$\bullet$ $E:=H_{r}^{1}(\R^{2})\times H_{r}^{1}(\R^{2})$ with norm
\begin{align*}
\|(u,v)\|_{E}^{2}=\|u\|_{H_{r}^{1}(\R^{2})}^{2}+\|v\|_{H_{r}^{1}(\R^{2})}^{2}.
\end{align*}
$\bullet$ $L^{2p}(\R^{2})\times L^{2p}(\R^{2})$ for $p>1$ with the norm
\begin{align*}
\|(u,v)\|_{2p}^{2p}=\|u\|_{2p}^{2p}+\|u\|_{2p}^{2p}.
\end{align*}
$\bullet$ $\mathbb{L}^{2p}_{loc}(\R^{2}) :=L^{2p}_{loc}(\R^{2})\times L^{2p}_{loc}(\R^{2})$;

\bl
\label{Lemma 2.1}\cite{J-H-J}
{\it Suppose that a sequence ${u_{n}}$ converges weakly to a function $u$ in $H_{r}^{1}(\R^{2})$ as $n\rightarrow\infty$. Then for each $\varphi\in H_{r}^{1}(\R^{2})$, $B'(u_{n})$, $B'(u_{n})\varphi$ and $B'(u_{n})u_{n}$ converges up to a subsequence to $B(u)$, $B'(u)\varphi$ and $B'(u)u$, respectively, as $n\rightarrow\infty$.}
\el

\bl
\label{Lemma 2.2}\cite{J-H-J} {\it For $u\in H_{r}^{1}(\R^{2})$, the following inequality holds
\begin{align*}
\int_{\R^{2}}|u|^{4}dx\leq4\Big(\int_{\R^{2}}|\nabla u|^{2}dx\Big)^{\frac{1}{2}}\Big(\int_{\R^{2}}\frac{u^{2}}{|x|^{2}}\Big(\int_{0}^{|x|}\frac{s}{2}u^{2}(s)ds\Big)^{2}dx\Big)^{\frac{1}{2}}.
\end{align*}
Furthermore, the equality is attained by a continuum of functions
\begin{align*}
\Big\{u_{l}=\frac{\sqrt{8}l}{1+|lx|^{2}}\in H_{r}^{1}(\R^{2})\Big|\,\, l\in(0,\infty)\Big\},
\end{align*}
\begin{align*}
\frac{1}{4}\int_{\R^{2}}|u_{l}|^{4}dx=\int_{\R^{2}}|\nabla u_{l}|^{2}dx=\int_{\R^{2}}\frac{u_{l}^{2}}{|x|^{2}}\Big(\int_{0}^{|x|}\frac{s}{2}u_{l}^{2}(s)ds\Big)^{2}dx=\frac{16\pi l^{2}}{3}.
\end{align*}
}
\el

\bl
\label{Lemma 2.3}
{\it If $(u,v)$ is a solution of (\ref{question1}) then it satisfies the Pohozaev identity
\begin{equation}
\label{3}
\|u\|_{2}^{2}+\omega\|v\|_{2}^{2}+2\int_{\R^{2}}\bigg(\frac{h^{2}(|x|)}{|x|^{2}}u^{2}+\frac{g^{2}(|x|)}{|x|^{2}}v^{2}\bigg)dx=\frac{1}{p}\Big(\|(u,v)\|_{2p}^{2p}+2b\|uv\|_{p}^{p}\Big).
\end{equation}
}
\el

\noindent{\it Proof.}
We adopt some idea in \cite{J-H-J}. Assume that $(u,v)\in E$ is a weak solution for problem (\ref{question1}). Similar to \cite{AA} and \cite{J-H-J}, we know that $\int_{|x|}^{\infty}\frac{h(s)}{s}u^{2}(s)ds,\int_{|x|}^{\infty}\frac{g(s)}{s}v^{2}(s)ds, \frac{h^{2}(|x|)}{|x|^{2}},\frac{g^{2}(|x|)}{|x|^{2}}\in
L^{\infty}(\R^{2})$.
Thus, the standard elliptic estimates \cite{GT} imply that $u,v\in C_{loc}^{1,\gamma}(\R^{2})$ for some $\gamma>0$. Then we obtain $\int_{|x|}^{\infty}\frac{h(s)}{s}u^{2}(s)ds, \int_{|x|}^{\infty}\frac{g(s)}{s}v^{2}(s)ds, \frac{h^{2}(|x|)}{|x|^{2}},\frac{g^{2}(|x|)}{|x|^{2}}\in C(\R^{2})$. Since $u,v\in H_{r}^{1}(\R^{2})$, we deduce that $u,v\in C^{2}(\R^{2})$. Then, multiplying the first equation in \eqref{question1} by $x\cdot\nabla u$ and integrating by parts on a ball $B_{R}=\{x\in\R^{2}:|x|<R\}$, then
\begin{align*}
&\int_{B_{R}}\Delta u(\nabla u\cdot x)dx=\frac{R}{2}\int_{\partial B_{R}}|\nabla u|^{2}dS_{x},\\
&\int_{B_{R}}u(\nabla u\cdot x)dx=-\int_{B_{R}}u^{2}dx+o_{R}(1),\\
&\int_{B_{R}}|u|^{2p-2}u(\nabla u\cdot x)dx=-\frac{1}{p}\int_{B_{R}}u^{2p}dx+o_{R}(1),\\
&\int_{B_{R}}\bigg(\int_{|x|}^{\infty}\frac{h(s)}{s}u^{2}(s)ds\bigg)u(\nabla u\cdot x)dx+\int_{B_{R}}\frac{h^{2}(|x|)}{|x|^{2}}u(\nabla u\cdot x)dx=\pi h^{2}(R)u^{2}(R)\\
&\ \ \ \ \ \ \
+\pi\bigg(\int_{R}^{+\infty}\frac{h(s)}{s}u^{2}(s)ds\bigg)u^{2}(R)R^{2}-2\int_{\R^{2}}\frac{h^{2}(|x|)}{|x|^{2}}u^{2}dx+o_{R}(1).
\end{align*}
Thus, we deduce that
\begin{equation}\label{4}
\begin{split}
&\frac{R}{2}\int_{\partial B_{R}}|\nabla u|^{2}dS_{x}+\int_{B_{R}}u^{2}dx-\pi h^{2}(R)u^{2}(R)-\pi\bigg(\int_{R}^{+\infty}\frac{h(s)}{s}u^{2}(s)ds\bigg)u^{2}(R)R^{2}\\ &\ \ \ \ \ +2\int_{\R^{2}}\frac{h^{2}(|x|)}{|x|^{2}}u^{2}dx=\frac{1}{p}\int_{B_{R}}u^{2p}dx-b\int_{B_{R}}|v|^{p}|u|^{p-2}ux\cdot\nabla udx.
\end{split}
\end{equation}
In a similar way,
\begin{equation}\label{5}
\begin{split}
&\frac{R}{2}\int_{\partial B_{R}}|\nabla v|^{2}dS_{x}+\int_{B_{R}}\omega v^{2}dx-\pi g^{2}(R)v^{2}(R)-\pi\bigg(\int_{R}^{+\infty}\frac{g(s)}{s}v^{2}(s)ds\bigg)v^{2}(R)R^{2}\\
&\ \ \ \ \ +2\int_{\R^{2}}\frac{g^{2}(|x|)}{|x|^{2}}v^{2}dx=\frac{1}{p}\int_{B_{R}}v^{2p}dx-b\int_{B_{R}}|u|^{p}|v|^{p-2}vx\cdot\nabla vdx,
\end{split}
\end{equation}
can be obtained. Then, summing up (\ref{4}) and (\ref{5}), we obtain
\begin{align*}
&\int_{B_{R}}(u^{2}+\omega v^{2})dx+2\int_{\R^{2}}\bigg(\frac{h^{2}(|x|)}{|x|^{2}}u^{2}+\frac{g^{2}(|x|)}{|x|^{2}}v^{2}\bigg)dx-\frac{1}{p}\int_{B_{R}}(u^{2p}+v^{2p})dx-\frac{2b}{p}\int_{B_{R}}|uv|^{p}dx\\
&\ \ \ =\pi\bigg(\int_{R}^{+\infty}\frac{h(s)}{s}u^{2}(s)ds\bigg)u^{2}(R)R^{2}+\pi\bigg(\int_{R}^{+\infty}\frac{g(s)}{s}v^{2}(s)ds\bigg)v^{2}(R)R^{2}+\pi h^{2}(R)u^{2}(R)\\
&\ \ \ \ \ \ \ +\pi g^{2}(R)v^{2}(R)-\frac{R}{2}\int_{\partial B_{R}}(|\nabla u|^{2}+|\nabla v|^{2})dS_{x}.
\end{align*}
Arguing as in \cite[Proposition 2.3.]{J-H-J}, there exists a suitable sequence $R_{n}\rightarrow\infty$ on which the right hand side above tends to zero. Passing to the limit we get the identity. This completes the proof.
\qed

\

\noindent{\it Proof of Theorem \ref{Theorem 1.2.}.}
Let $(u,v)$ be a solution of \eqref{question1}. By Lemma \ref{Lemma 2.2}, we obtain
\begin{align*}
0&=\int_{\R^{2}}(|\nabla u|^{2}+u^{2}+|\nabla v|^{2}+\omega v^{2})dx+3\Big(B(u)+B(v)\Big)-\int_{\R^{2}}(u^{2p}+v^{2p}+2b|uv|^{p})dx\\
&\geq\int_{\R^{2}} (u^{2}+\frac{1}{2}u^{4}-(1+b)u^{2p})dx+\int_{\R^{2}}(\omega v^{2}+\frac{1}{2}v^{4}-(1+b)v^{2p})dx.
\end{align*}
Denote
\begin{align*}
&f_{1}(t)=t^{2}+\frac{1}{2}t^{4}-(1+b)t^{2p},\\
&f_{1}'(t)=2t+2t^{3}-2p(1+b)t^{2p-1}=2t\Big(1+t^{2}-2p(1+b)t^{2p-2}\Big),
\end{align*}
There exists $\tilde{b}>0$ small enough such that
\begin{align*}
1+\Big(p(1+\tilde{b})\Big)^{\frac{1}{2-p}}\Big(p-1\Big)^{\frac{p-1}{2-p}}(p-2)=0.
\end{align*}
Then, we have
\begin{align*}
f_{1}(t)\geq0,\ \ \ t\in\mathbb{R},\ \ \ b\in(0,\tilde{b}).
\end{align*}
There exists a $\tilde{\omega}>0$ such that the function $t\mapsto \omega t^{2}+\frac{1}{2}t^{4}-(1+b)t^{2p}$ is nonnegative and strictly increases as $\omega>\tilde{\omega},b\in(0,\tilde{b})$. Hence, $(u,v)$ must be identically $(0,0)$. This completes the proof.
\qed

\vs
\section{Proof of $(ii)$ and $(iii)$ of Theorem \ref{Theorem 1.1.}. }
\renewcommand{\theequation}{3.\arabic{equation}}
Consider the following problem
\begin{align}\label{singal}
-\Delta u+\omega u+\big(\int_{|x|}^{\infty}\frac{h(s)}{s}u^{2}(s)ds+\frac{h^{2}(|x|)}{|x|^{2}}\big)u=|u|^{2p-2}u,\ x\in\R^{2},
\end{align}
By \cite{J-H-J}, when $p\in(3,\infty)$, problem \eqref{singal} admits a positive ground state solution $u_{\omega}$. To be more precise, define the associated energy functional by
\begin{align*}
J_{\omega}=\frac{1}{2}\int_{\R^{2}}(|\nabla u|^{2}+\omega u^{2})dx+\frac{1}{2}B(u)-\frac{1}{2p}\int_{\R^{2}}u^{2p}dx.
\end{align*}
Denote the ground state level by
\begin{align*}
E_{\omega}:=\min\{J_{\omega}(u):u\in H_{r}^{1}(\R^{2})\setminus\{0\},J'_{\omega}(u)=0\},
\end{align*}
Moreover,
\begin{align*}
E_{\omega}:=\inf_{u\in\mathcal{N}_{\omega}}J_{\omega}(u),
\end{align*}
where
\begin{align*}
\mathcal{N}_{\omega}:=\Bigg\{u\in H_{r}^{1}(\R^{2})\setminus\{0\}:\int_{\R^{2}}(|\nabla u|^{2}+\omega u^{2})dx+3B(u)=\int_{\R^{2}}u^{2p}dx\Bigg\}.
\end{align*}

Define the Nehari manifold of problem (\ref{question1}) by
\bg
\N:=\left\{(u,v)\in E\setminus\left\{(0,0)\right\}|\langle I'(u,v),(u,v)\rangle=0\right\}.
\end{align*}
The corresponding groundstate energy is described as
\begin{align*}
c_{\N}:=\mathop{\inf}\limits_{(u,v)\in\N}I(u,v).
\end{align*}\
\bl
\label{Lemma 3.1}{\it(Mountain-Pass geometry)} {\it Assume $b>0$, then the functional $I$ satisfies the following conditions\\
$(i)$ There exists a positive constant $r>0$ such that $I(u,v)>0$ for $\|(u,v)\|_E=r$;\\
$(ii)$ There exists $(e_{1},e_{2})\in E$ with $\|(e_{1},e_{2})\|_{E}>r$ such that $I(e_{1},e_{2})<0$.\
}
\el

\noindent{\it Proof.} Since
\begin{align*}
2b\int|uv|^{p}dx\leq b(\|u\|_{2p}^{2p}+\|v\|_{2p}^{2p}),
\end{align*}
and by the Sobolev embedding theorem, there exists a positive constant $C$ such that
\begin{align*}
I(u,v)&=\frac{1}{2}\|(u,v)\|_{E}^{2}+\frac{1}{2}\Big(B(u)+B(v)\Big)-F(u,v);\\
&\geq\frac{1}{2}\|(u,v)\|_{E}^{2}-C\|(u,v)\|_{E}^{2p}.
\end{align*}
Hence, there exists $r>0$ such that
\begin{align*}
\mathop{\inf}\limits_{\|(u,v)\|_{E}=r}I(u,v)>0.
\end{align*}
For any $(u,v)\in E\setminus\{(0,0)\}$ and $t>0$,
\begin{align*}
I(tu,tv)=\frac{t^{2}}{2}\|(u,v)\|_{E}^{2}+\frac{t^{6}}{2}\Big(B(u)+B(v)\Big)-\frac{t^{2p}}{2p}F(u,v),
\end{align*}
which implies that $I(tu,tv)\rightarrow-\infty$ as $t\rightarrow+\infty$. This completes the proof.
\qed

That is $I$ satisfies the geometric conditions of the Mountain-Pass theorem. Define
\begin{align*}
c=\mathop{\inf}\limits_{\gamma\in\Gamma}\mathop{\max}\limits_{t\in[0,1]}I(\gamma(t)),
\end{align*}
where $\Gamma=\left\{\gamma\in C([0,1],E):I(\gamma(0))=0, I(\gamma(1))<0\right\}$.
\

\bl
\label{Lemma 3.2}{\it For every $(u,v)\in E\setminus\{(0,0)\}$ there exists a unique $\overline{t}_{uv}>0$ such that $\overline{t}_{uv}(u,v)\in\N$. The maximum of $I(tu,tv)$ for $t\geq0$ is achieved at $t=\overline{t}_{uv}$.}
\el

\noindent{\it Proof.} $\forall(u,v)\in E\setminus\{(0,0)\}$ and $t>0$, let
\begin{align*}
g(t):=I(tu,tv)=\frac{t^{2}}{2}\|(u,v)\|_{E}^{2}+\frac{t^{6}}{2}\Big(B(u)+B(v)\Big)-\frac{t^{2p}}{2p}F(u,v).
\end{align*}
By Lemma \ref{Lemma 3.1}, we deduce that there exists $\overline{t}=\overline{t}_{uv}>0$ such that
\begin{align*}
g(\overline{t})=\mathop{\max}\limits_{t>0}g(t).
\end{align*}
Moreover
\begin{align*}
g'(t)&=t\|(u,v)\|_{E}^{2}+3t^{5}\Big(B(u)+B(v)\Big)-t^{2p-1}F(u,v)\\
&=t^{2p-1}\Big(\frac{1}{t^{2p-2}}\|(u,v)\|_{E}^{2}+\frac{3}{t^{2p-6}}\big(B(u)+B(v)\big)-F(u,v)\Big).
\end{align*}
As $p>3$, thus $g'(t)$ is strictly decreasing, the point $t=\overline{t}_{uv}$ is the unique value of $t>0$ at which $\overline{t}_{uv}(u,v)\in\N$. This completes the proof.
\qed
\

\bl\label{Lemma 3.3}
{\it $c_{\N}=c$}.
\el

\noindent{\it Proof.} Define
\begin{align*}
c_{1}:=\mathop{\inf}\limits_{(u,v)\in E\setminus\{(0,0)\}}\mathop{\max}\limits_{t\geq0}I(tu,tv).
\end{align*}
By Lemma \ref{Lemma 3.2}, we can obtain $c_{\N}=c_{1}$. Since $I(tu,tv)<0$ for any $t$ large, it follows that $c\leq c_{1}$. Since every $\gamma\in\Gamma$ intersects $\N$, so that $c\geq c_{\N}$. This completes the proof.
\qed

\

Set
\begin{align*}
&b_{1}=\Bigg(\frac{(p-1)3^{\frac{p}{p-1}}\|u_{1}\|_{2p}^{2p}}{pE_{\omega}}\Bigg)^{p-1},\ \ \ b_{2}=\Bigg(\frac{(p-1)3^{\frac{p}{p-1}}\|u_{\omega}\|_{2p}^{2p}}{pE_{1}}\Bigg)^{p-1},\\
&b_{3}=\Bigg(\frac{(p-1)3^{\frac{3}{p-3}}\|u_{1}\|_{2p}^{2p}}{pE_{\omega}}\Bigg)^{p-1},\ \ \ b_{4}=\Bigg(\frac{(p-1)3^{\frac{3}{p-3}}\|u_{\omega}\|_{2p}^{2p}}{pE_{1}}\Bigg)^{p-1}.
\end{align*}
\

\bl
{Assume that one of the following conditions holds\\
$(i)$ $p\in(3,3+\sqrt{6})$ and $b>\max\{b_{1},b_{2}\}$,\\
$(ii)$ $p\in[3+\sqrt{6},\infty)$ and $b>\max\{b_{3},b_{4}\},$\\
then $c_{\mathcal{N}}<\min\{E_{1},E_{\omega}\}$, where $E_{1}$ is a groundstate energy of problem \eqref{singal} with $\omega=1$.}
\el

\noindent{\it Proof.}
Denote
\begin{align*}
&a(u,v):=\|\nabla u\|_{2}^{2}+\|\nabla v\|_{2}^{2},\\
&b(u,v):=\|u\|_{2}^{2}+\omega\|v\|_{2}^{2},\\
&c(u,v):=B(u)+B(v).
\end{align*}
For fixed $(u,v)\in E\setminus\{(0,0)\}$, we have
\begin{align*}
c_{\N}&\leq\mathop{\max}\limits_{t\geq0}I(tu,tv)\\
&\leq\mathop{\max}\limits_{t\geq0}\Bigg\{\frac{t^{2}}{2}\big(a(u,v)+b(u,v)\big)+\frac{t^{6}}{2}c(u,v)-\frac{t^{2p}}{2p}F(u,v)\Bigg\}\\
&\leq\mathop{\max}\limits_{t\geq0}\Bigg\{\frac{t^{2}}{2}\big(a(u,v)+b(u,v)\big)-\frac{t^{2p}}{4p}F(u,v)\Bigg\}+\mathop{\max}\limits_{t\geq0}\Bigg\{\frac{t^{6}}{2}c(u,v)-\frac{t^{2p}}{4p}F(u,v)\Bigg\}\\
&\leq\frac{p-1}{2p}\Bigg(\frac{2(a(u,v)+b(u,v))^{p}}{F}\Bigg)^{\frac{1}{p-1}}+\frac{p-3}{2p}\Bigg(\frac{6c(u,v)^{\frac{p}{3}}}{F}\Bigg)^{\frac{3}{p-3}}.
\end{align*}
Assume $\omega\leq1$, $E_{1}\geq E_{\omega}$. Choosing $(u_{1},u_{1})$ such that
\begin{align*}
c_{\N}\leq \mathop{\max}\limits_{t\geq0}I(tu_{1},tu_{1})<E_{\omega},
\end{align*}
where $u_{1}$ is a positive ground state solution of \eqref{singal} with $\omega=1$. Then
\begin{align*}
c_{\N}&\leq\mathop{\max}\limits_{t\geq0}I(tu_{1},tu_{1})\\
&\leq\frac{p-1}{2p}\Bigg(\frac{2(2\|u_{1}\|^{2})^{p}}{2(1+b)\|u_{1}\|_{2p}^{2p}}\Bigg)^{\frac{1}{p-1}}+\frac{p-3}{2p}\Bigg(\frac{6(2B(u_{1}))^{\frac{p}{3}}}{2(1+b)\|u_{1}\|_{2p}^{2p}}\Bigg)^{\frac{3}{p-3}}\\
&\leq\max\{2^{\frac{p}{p-1}},3^{\frac{3}{p-3}}\}\frac{p-1}{p}\|u_{1}\|_{2p}^{2p}\Big(\frac{1}{1+b}\Big)^{\frac{1}{p-1}}\\
&<\max\{2^{\frac{p}{p-1}},3^{\frac{3}{p-3}}\}\frac{p-1}{p}\|u_{1}\|_{2p}^{2p}\Big(\frac{1}{b}\Big)^{\frac{1}{p-1}}
\end{align*}
Due to
\begin{align*}
\max\{2^{\frac{p}{p-1}},3^{\frac{3}{p-3}}\}=
\bcs
3^{\frac{p}{p-1}},\ \ \ p\in[3+\sqrt{6},\infty),\\
3^{\frac{3}{p-3}},\ \ \ p\in(3,3+\sqrt{6}),
\ecs
\end{align*}
we deduce that
\begin{align*}
b>
\bcs
\Bigg(\frac{(p-1)3^{\frac{p}{p-1}}\|u_{1}\|_{2p}^{2p}}{pE_{\omega}}\Bigg)^{p-1},\ \ \ p\in[3+\sqrt{6},\infty),\\
\Bigg(\frac{(p-1)3^{\frac{3}{p-3}}\|u_{1}\|_{2p}^{2p}}{pE_{\omega}}\Bigg)^{p-1},\ \ \ p\in(3,3+\sqrt{6}).
\ecs
\end{align*}
In a similar fashion, if $\omega\geq1$, it follows that
\begin{align*}
b>
\bcs
\Bigg(\frac{(p-1)3^{\frac{p}{p-1}}\|u_{\omega}\|_{2p}^{2p}}{pE_{1}}\Bigg)^{p-1},\ \ \ p\in[3+\sqrt{6},\infty),\\
\Bigg(\frac{(p-1)3^{\frac{3}{p-3}}\|u_{\omega}\|_{2p}^{2p}}{pE_{1}}\Bigg)^{p-1},\ \ \ p\in(3,3+\sqrt{6}).
\ecs
\end{align*}
In conclusion, $c_{\N}<\min\{E_{1},E_{\omega}\}$ provided
\begin{align*}
b>
\bcs
\max\{b_{1},b_{2}\}\ \ \  p\in[3+\sqrt{6},\infty),\\
\max\{b_{3},b_{4}\}\ \ \  p\in(3,3+\sqrt{6}).
\ecs
\end{align*}
This completes the proof.
\qed

\

{\bf \it Proof of $(ii)$ and $(iii)$ of Theorem \ref{Theorem 1.1.}.} We divide the proof in several steps.

{\bf Step 1.}
We show that the existence of ground state solutions of problem (\ref{question1}). By the Ekeland variational principle \cite{Willem} that there exists a sequence $\{(u_{n},v_{n})\}\in E$ such that
\begin{align*}
I(u_{n},v_{n})\rightarrow c,\ \ \ I'(u_{n},v_{n})\rightarrow0\ \mbox{in}\ E'.
\end{align*}
By computing $I(u_{n},v_{n})-\frac{1}{6}\langle I'(u_{n},v_{n}),(u_{n},v_{n})\rangle$, it is easy to get that $\{(u_{n},v_{n})\}$ is bounded in $E$. Then there exists $(u,v)\in E$ such that, up to a subsequence
\begin{align*}
&(u_{n},v_{n})\rightarrow(u,v)\ \mbox{weakly\ in}\ E,\\
&(u_{n},v_{n})\rightarrow(u,v)\ \mbox{strongly\ in}\ L^{2p}(\R^{2})\times L^{2p}(\R^{2})\ \mbox{for}\
 p\in(3,+\infty),\\
&(u_{n},v_{n})\rightarrow(u,v)\ \mbox{a.e.}\ \mbox{in}\ \R^{2}.
\end{align*}
Then $I'(u,v)=0$ and $I(u,v)\leq c$. If $(u,v)\neq(0,0)$, it follows that $(u,v)\in\N$, $I(u,v)\geq c$. Then, $I(u,v)=c$. Now, it remains to prove $(u,v)\neq(0,0)$. Assume by the contrary that $(u,v)=(0,0)$, then
\begin{align*}
c+o_{n}(1)\|(u_{n},v_{n})\|_{E}&=I((u_{n},v_{n}))-\frac{1}{2}\langle I'(u_{n},v_{n}),(u_{n},v_{n})\rangle\\
&=-\big(B(u+B(v)\big)+\Big(\frac{1}{2}-\frac{1}{2p}\Big)F(u,v)+o_{n}(1)\\
&=o_{n}(1),
\end{align*}
which is a contradiction. Thus $(u,v)$ is a ground state solution of system (\ref{question1}).

{\bf Step 2.}
We show that $u\neq0$ and $v\neq0$. Without loss of generality, we assume that $u=0$ and $v\not\equiv0$. Multiplying the first equation in \eqref{question1} by $u_{n}$, we obtain
\begin{align*}
\lim_{n\rightarrow\infty}\int_{\R^{2}}(|\nabla u_{n}|^{2}+u_{n}^{2})dx=0,
\end{align*}
Multiplying the second equation in \eqref{question1} by $v_{n}$,
\begin{align}\label{10}
\int_{\R^{2}}(|\nabla v_{n}|^{2}+v_{n}^{2})dx+3B(v_{n})=\int_{\R^{2}}v_{n}^{2p}dx+o_{n}(1).
\end{align}
Therefore, there exists $t_{n}$ such that
\begin{align}\label{15}
\frac{1}{t_{n}^{2}}\int_{\R^{2}}(|\nabla v_{n}|^{2}+v_{n}^{2})dx+3B(v_{n})=t_{n}^{2p-6}\int_{\R^{2}}v_{n}^{2p}dx.
\end{align}
Combining \eqref{10} and \eqref{15}, it follows that $t_{n}\rightarrow1$, as $n\rightarrow\infty$. Hence
\begin{align*}
\lim_{n\rightarrow\infty}I(u_{n},v_{n})\rightarrow I(0,v)\geq E_{\omega},
\end{align*}
which contradicts the fact that $c_{\N}<E_{\omega}$. Similarly, if $v=0$ and $u\not\equiv0$, we can obtain $I(u_{n},v_{n})\rightarrow I(u,0)$, as $n\rightarrow\infty$, which contradicts the fact that $c_{\N}<E_{1}$.

Therefore, $u\not\equiv0$ and $v\not\equiv0$, $(u,v)$ is a nontrivial ground state solution of \eqref{question1}. In fact, since $(|u|,|v|)\in \N$ and $c_{\N}=I(|u|,|v|)$, we conclude that $(|u|,|v|)$ is a nonnegative solution of (\ref{question1}). Using the strong maximum principle, we infer that $|u|,|v|>0$. Thus $(|u|,|v|)$ is a positive least energy solution of (\ref{question1}). This completes the proof.
\qed

\vs
\section{Proof of $(i)$ of Theorem \ref{Theorem 1.1.}.}
\renewcommand{\theequation}{4.\arabic{equation}}
Given $(u,v)\in E\setminus\{(0,0)\}$, consider the path
\begin{align*}
\gamma_{u,v}(t):=(t^{\alpha}u(t\cdot),t^{\alpha}v(t\cdot)),\ \ t\geq0,
\end{align*}
where $\alpha>1$ such that$\frac{1}{p-1}<\alpha<\frac{1}{3-p}$ for $p\in(2,3)$ and $\alpha>1$ for $p=3$.
Then
\begin{align*}
I(\gamma_{u,v}(t))=&\frac{t^{2\alpha}}{2}\big(\|\nabla u\|_{2}^{2}+\|\nabla v\|_{2}^{2}\big)+\frac{t^{2(\alpha-1)}}{2}\big(\|u\|_{2}^{2}+\omega\|v\|_{2}^{2}\big)\\
&+\frac{t^{6\alpha-4}}{2}\Big(B(u)+B(v)\Big)-\frac{t^{2p\alpha-2}}{2p}F(u,v).
\end{align*}
By differentiating both sides with respect to $t$ at
$1$, we obtain the following constraint
\begin{align*}
J(u,v)=&\alpha\big(\|\nabla u\|_{2}^{2}+\|\nabla v\|_{2}^{2}\big)+(\alpha-1)\big(\|u\|_{2}^{2}+\omega\|v\|_{2}^{2}\big)\\
&+(3\alpha-2)\Big(B(u)+B(v)\Big)-\frac{p\alpha-1}{p}F(u,v).
\end{align*}
Define a constraint manifold of Pohozaev-Nehari mype
\begin{align*}
\M_{b}:=\left\{(u,v)\in E\setminus\left\{(0,0)\right\}|J(u,v)=0\right\}.
\end{align*}
The corresponding groundstate energy is described as
\begin{align*}
c_{b}:=\mathop{\inf}\limits_{(u,v)\in\M_{b}}I(u,v).
\end{align*}

By \cite{J-H-J}, when $p\in(2,3]$, problem \eqref{singal} admits a positive ground state solution $\widetilde{u}_{\omega}$. To be more precise, denote the ground state level by
\begin{align*}
\widetilde{E}_{\omega}:=\inf_{u\in\mathcal{M}_{\omega}}J_{\omega}(u)
\end{align*}
where
\begin{align*}
\mathcal{M}_{\omega}:=\Bigg\{u\in H_{r}^{1}(\R^{2})\setminus\{0\}:\int_{\R^{2}}(\alpha|\nabla u|^{2}+(\alpha-1)\omega u^{2})dx+(3\alpha-2)B(u)=\frac{p\alpha-1}{p}\int_{\R^{2}}u^{2p}dx\Bigg\}.
\end{align*}
Define the set of ground state solutions of \eqref{singal} by
\begin{align*}
S_{\omega}:=\{u\in H_{r}^{1}(\R^{2})\setminus\{0\}:J'_{\omega}(u)=0,J_{\omega}(u)=E_{\omega}\}.
\end{align*}
$S_{\omega}$ is nonempty.

\bl\label{Lemma 4.1}
{The set $S_{\omega}$ is compact in $H_{r}^{1}(\R^{2})$. More precisely, any sequence $\{u_{n}\}\subset S_{\omega}$ has a subsequence $\{u_{jn}\}\subset S_{\omega}$ and $u\in S_{\omega}$ such that $u_{jn}\rightarrow u$ strongly in $H_{r}^{1}(\R^{2})$ as $n\rightarrow\infty$.}
\el

\noindent{\it Proof.}
For any  $\{u_{n}\}\in S_{\omega}$,
\begin{align*}
E_{\omega}&=\frac{1}{2}\int_{\R^{2}}(|\nabla u_{n}|^{2}+\omega u_{n}^{2})dx+\frac{1}{2}B(u_{n})-\frac{1}{2p}\int_{\R^{2}}u_{n}^{2p}dx\\
&=\Big(\frac{1}{2}-\frac{\alpha}{2(p\alpha-1)}\Big)\int_{\R^{2}}|\nabla u_{n}|^{2}dx+\Big(\frac{1}{2}-\frac{\alpha-1}{2(p\alpha-1)}\Big)\omega\int_{\R^{2}}u_{n}^{2}dx+\Big(\frac{1}{2}-\frac{3\alpha-2}{2(p\alpha-1)}\Big)B(u_{n})\\
&\geq\Big(\frac{1}{2}-\frac{\alpha}{2(p\alpha-1)}\Big)\int_{\R^{2}}(|\nabla u_{n}|^{2}+u_{n}^{2})dx.
\end{align*}
Therefore $\{u_{n}\}$ is bounded on $S_{\omega}$. There exists $u\in H_{r}^{1}(\R^{2})$, such that, up to a subsequence, $u_{n}\rightarrow u$ weakly in $H_{r}^{1}(\R^{2})$, strongly in $L^{p}(\R^{2})$ for $p\in(2,3]$ as $n\rightarrow\infty$. Moreover, due to $u_{n}\in\mathcal{M}_{\omega}$ for any $n$, we can obtain
\begin{align}\label{11}
\liminf_{n\rightarrow\infty}\|u_{n}\|_{H_{r}^{1}(\R^{2})}>0.
\end{align}
First we prove that $u\not\equiv0$. Indeed, if $u=0$, by Lemma \ref{Lemma 2.1}, we obtain $u_{n}\rightarrow0$ strongly in $H_{r}^{1}(\R^{2})$, as $n\rightarrow\infty$, which contradicts \eqref{11}. Thus $u\not\equiv0$. Next, we show that $u\in S_{\omega}$. Since $u_{n}$ satisfies problem \eqref{singal}, we get that $u$ is a solution of problem \eqref{singal}. By the semicontinuity of the norms,
\begin{align*}
\widetilde{E}_{\omega}\leq J_{\omega}(u)\leq\liminf_{n\rightarrow\infty}J_{\omega}(u_{n})=\widetilde{E}_{\omega}.
\end{align*}
Therefore $u\in S_{\omega}$ and $u_{n}\rightarrow u$ strongly in $H_{r}^{1}(\R^{2})$, as $n\rightarrow\infty$. That is, $S_{\omega}$ is compact in $H_{r}^{1}(\R^{2})$.
\qed

\bl\label{Lemma 3.4}
{\it For given positive constants $a$, $b$, $c$, $d$, a function $f(t)=at^{2\alpha}+bt^{2(\alpha-1)}+ct^{6\alpha-4}-dt^{2p\alpha-2}$ has exactly one critical point on $(0,+\infty)$, the maximum.}
\el

\noindent{\it Proof.} The proof is similar to \cite{J-H-J}, we omit it here.
\qed

\bl\label{Lemma 3.5}
{\it For any $(u,v)\in E\setminus\{(0,0)\}$, there exists a unique $t_{u,v}>0$ such that $\gamma_{u,v}(t_{u,v})\in\M_{b}$ and
\begin{align*}
c_{b}=\mathop{\inf}\limits_{(u,v)\in E\setminus\{(0,0)\}}\mathop{\max}\limits_{t>0}I(\gamma_{u,v}(t)).
\end{align*}
}
\el

\noindent{\it Proof.} The proof is standard, we omit it here.
\qed

\bl\label{Lemma 3.6}
{If $b>0$, then $c_{b}<\widetilde{E}_{1}+\widetilde{E}_{\omega}$, where $\widetilde{E}_{1}$ is a groundstate energy of problem \eqref{singal} with $\omega=1$.
}
\el

\noindent{\it Proof.} Let $\widetilde{u}_{1}$ and $\widetilde{u}_{\omega}$ be a positive ground state solution associated to the level $\widetilde{E}_{1}$ and $\widetilde{E}_{\omega}$ respectively. Denote $f(t)=I(t^{\alpha}\widetilde{u}_{1}(t\cdot),t^{\alpha}\widetilde{u}_{\omega}(t\cdot))$, that is,\
\begin{align*}
f(t)=&\frac{t^{2\alpha}}{2}\big(\|\nabla \widetilde{u}_{1}\|_{2}^{2}+\|\nabla \widetilde{u}_{\omega}\|_{2}^{2}\big)+\frac{t^{2(\alpha-1)}}{2}\big(\|\widetilde{u}_{1}\|_{2}^{2}+\omega\|\widetilde{u}_{\omega}\|_{2}^{2}\big)\\
&+\frac{t^{6\alpha-4}}{2}\Big(B(\widetilde{u}_{1})+B(\widetilde{u}_{\omega})\Big)-\frac{t^{2p\alpha-2}}{2p}F(\widetilde{u}_{1},\widetilde{u}_{\omega}).
\end{align*}
According to Lemma \ref{Lemma 3.4}, there exists unique $t_{0}>0$ such that $(t_{0}^{\alpha}\widetilde{u}_{1}(t_{0}\cdot),t_{0}^{\alpha}\widetilde{u}_{\omega}(t_{0}\cdot))\in\M_{b}$. Hence
\begin{align*}
c_{b}&\leq I(t_{0}^{\alpha}\widetilde{u}_{1}(t_{0}\cdot),t_{0}^{\alpha}\widetilde{u}_{\omega}(t_{0}\cdot))=I(t_{0}\widetilde{u}_{1}(t_{0}\cdot),0)+I(0,t_{0}\widetilde{u}_{\omega}(t_{0}\cdot))-2bt_{0}^{2p}\int_{\R^{2}}|\widetilde{u}_{1}(t_{0}x)\widetilde{u}_{\omega}(t_{0}x)|^{p}dx\\
&<I(t_{0}\widetilde{u}_{1}(t_{0}\cdot),0)+I(0,t_{0}\widetilde{u}_{\omega}(t_{0}\cdot))\\
&=\widetilde{E}_{1}+\widetilde{E}_{\omega}.
\end{align*}
This completes the proof.
\qed

\bl\label{Lemma 3.7}
{$\liminf_{b\rightarrow0}c_{b}>0$.
}
\el

\noindent{\it Proof.}
Suppose by contradiction the lemma does not hold. Then there exists $\{b_{k}\}$ such that $b_{k}\rightarrow0$ and $c_{b_{k}}\rightarrow0$, as $k\rightarrow\infty$. Moreover, there exists $\{(u_{k},v_{k})\}\in\mathcal{M}_{b_{k}}$ such that $I(u_{k},v_{k})\rightarrow0$, as $k\rightarrow\infty$, that is
\begin{align*}
o_{k}(1)&=\frac{1}{2}\int_{\R^{2}}(|\nabla u_{k}|^{2}+u_{k}^{2}+|\nabla v_{k}|^{2}+\omega v_{k}^{2})dx+\frac{1}{2}(B(u_{k})+B(v_{k}))-\frac{1}{2p}F(u_{k},v_{k})\\
&\geq\Big(\frac{1}{2}-\frac{\alpha}{2(p\alpha-1)}\Big)\|(u_{k},v_{k})\|_{E}^{2}.
\end{align*}
We deduce that $\|(u_{k},v_{k})\|_{E}\rightarrow0$, as $k\rightarrow\infty$. Since $\{(u_{k},v_{k})\}\in\mathcal{M}_{b_{k}}$,
\begin{align*}
(\al-1)\|(u_{k},v_{k})\|_{E}^{2}\leq C \|(u_{k},v_{k})\|_{E}^{2p}.
\end{align*}
Thus we have $\|(u_{k},v_{k})\|_{E}\geq(\frac{\al-1}{C})^{\frac{1}{2p-2}}$, which is a contradiction. This completes the proof.
\qed

\

Given $\delta>0$, let
\begin{align*}
(S_{\omega})^{\delta}:=\{u\in H_{r}^{1}(\R^{2})|u=\widetilde{u}+\overline{u},\widetilde{u}\in S_{\omega},\|\overline{u}\|_{H_{r}^{1}(\R^{2})}\leq\delta\},
\end{align*}
be the neighborhood of $S_{\omega}$ of radius $\delta$.

\bl\label{Lemma 3.8}
{For any $\delta>0$, there exists $b_{\delta}>0$ such that for any $b\in(0,b_{\delta})$, up to a subsequence, there exists ${(u_{n}^{b},v_{n}^{b})}\in\M_{b}$ satisfying
\begin{align}\label{6}
I(u_{n}^{b},v_{n}^{b})\rightarrow c_{b},\ \ I'(u_{n}^{b},v_{n}^{b})\rightarrow0, \mbox{as}\  n\rightarrow\infty,
\end{align}
and $\{u_{n}^{b}\}\in(S_{1})^{\delta}$ and $\{v_{n}^{b}\}\in(S_{\omega})^{\delta}$.
}
\el

{\it Proof.} We adopt some idea in \cite{DT}. Suppose by contradiction the lemma does not hold. Then for $\delta_{0}>0$, there exists $\{b_{k}\}\in\R^{+}$ such that $b_{k}\rightarrow0$, as $k\rightarrow\infty$, and for any $\{(u_{n}^{b_{k}},v_{n}^{b_{k}})\}\in\M_{b_{k}}$ satisfying (\ref{6}), there holds $\{u_{n}^{b_{k}}\}\in H_{r}^{1}(\R^{2})\setminus(S_{1})^{\delta_{0}}$ or $\{v_{n}^{b_{k}}\}\in H_{r}^{1}(\R^{2})\setminus(S_{\omega})^{\delta_{0}}$. For any $k$, there exists $n_{k}$ such that
\begin{align*}
|I(u_{n_{k}}^{b_{k}},v_{n_{k}}^{b_{k}})-c_{b_{k}}|\leq 1/k.
\end{align*}
Let $\widetilde{u}_{k}=u_{n_{k}}^{b_{k}}$ and $\widetilde{v}_{k}=v_{n_{k}}^{b_{k}}$. By Lemma \ref{Lemma 3.6}, we have
\begin{align}\label{20}
\limsup_{k\rightarrow\infty}I(\widetilde{u}_{k},\widetilde{v}_{k})\leq\limsup_{k\rightarrow\infty}c_{b_{k}}\leq \widetilde{E}_{1}+\widetilde{E}_{\omega}.
\end{align}
Since $\{(\widetilde{u}_{k},\widetilde{v}_{k})\}\in\mathcal{M}_{b_{k}}$,
\begin{align*}
c_{b_{k}}&=\frac{1}{2}\|(\widetilde{u}_{k},\widetilde{v}_{k})\|_{E}^{2}+\frac{1}{2}\Big(B(\widetilde{u}_{k})+B(\widetilde{v}_{k})\Big)-\frac{1}{2p}F(\widetilde{u}_{k},\widetilde{v}_{k})\\
&\geq\Big(\frac{1}{2}-\frac{\alpha}{2(p\alpha-1)}\Big)\|(\widetilde{u}_{k},\widetilde{v}_{k})\|_{E}^{2},
\end{align*}
we deduce that $\{(\widetilde{u}_{k},\widetilde{v}_{k})\}$ is bounded in $E$. Up to a subsequence, $\widetilde{u}_{k}\rightarrow u$ and $\widetilde{v}_{k}\rightarrow v$ weakly in $H_{r}^{1}(\R^{2})$, strongly in $L^{2p}(\R^{2})$ for $p\in(2,3]$, as $k\rightarrow\infty$. By Lemma \ref{Lemma 3.7}, we have
\begin{align*}
\liminf_{k\rightarrow\infty}\min\{\|\widetilde{u}_{k}\|_{H_{r}^{1}(\R^{2})},\|\widetilde{v}_{k}\|_{H_{r}^{1}(\R^{2})}\}>0.
\end{align*}
Noting that $b_{k}>0$ and
\begin{align*}
o_{k}(1)=\int_{\R^{2}}|\nabla\widetilde{u}_{k}|^{2}+\widetilde{u}_{k}^{2}+3B(\widetilde{u}_{k})-\int_{\R^{2}}\widetilde{u}_{k}^{2p},
\end{align*}
there exists $t_{k}$ such that
\begin{align*}
\frac{1}{t_{k}^{4\al-4}}\int_{\R^{2}}|\nabla\widetilde{u}_{k}|^{2}+\frac{1}{t_{k}^{4\al-2}}\int_{\R^{2}}\widetilde{u}_{k}^{2}+3B(\widetilde{u}_{k})=t_{k}^{2+2p\al-6\al}\int_{\R^{2}}\widetilde{u}_{k}^{2p},
\end{align*}
that is $t_{k}^{\alpha}\widetilde{u}_{k}(t_{k}\cdot)\in \mathcal{N}_{1}$. Similarly, there exists $s_{k}$ such that $s_{k}^{\alpha}\widetilde{v}_{k}(t_{k}\cdot)\in \mathcal{N}_{\omega}$.

{\bf Step 1.}
We claim that $t_{k}\rightarrow1$ and $s_{k}\rightarrow1$ as $k\rightarrow\infty$. We only give the proof of $t_{k}\rightarrow1$, as the second convergence being similar. We consider two cases:

{\bf Case I.}
$u\neq0$. If $\limsup_{k\rightarrow\infty}t_{k}>1$, then we can assume that $t_{k}>1$ for all $k$ we have
\begin{align*}
o_{k}(1)&=(t_{k}^{2+2p\al-6\al}-1)\int_{\R^{2}}\widetilde{u}_{k}^{2p}-(\frac{1}{t_{k}^{4\al-4}}-1)\int_{\R^{2}}|\nabla\widetilde{u}_{k}|^{2}-(\frac{1}{t_{k}^{4\al-2}}-1)\int_{\R^{2}}\widetilde{u}_{k}^{2}\\
&\geq(t_{k}^{2+2p\al-6\al}-1)\int_{\R^{2}}\widetilde{u}_{k}^{2p},
\end{align*}
which yields $t_{k}\rightarrow1$ as $k\rightarrow\infty$. This is a contradiction. So $\limsup_{k\rightarrow\infty}t_{k}\leq1$. Similarly, $\liminf_{k\rightarrow\infty}t_{k}\geq1$. Then $\lim_{k\rightarrow\infty}t_{k}=1$.

{\bf Case II.}
$u=0$. If $\limsup_{k\rightarrow\infty}t_{k}>1$, then we can assume that $t_{k}>1$ for all $k$ we have $\limsup_{k\rightarrow\infty}\|\widetilde{u}_{k}\|_{H_{r}^{1}(\R^{2})}=0$, which contradicts \eqref{20}. So $\limsup_{k\rightarrow\infty}t_{k}\leq1$. Similarly, $\liminf_{k\rightarrow\infty}t_{k}\geq1$. Then $\lim_{k\rightarrow\infty}t_{k}=1$.

{\bf Step 2.}
Let $\overline{u}_{k}=t_{k}^{\alpha}\widetilde{u}_{k}(t_{k}\cdot)$ and $\overline{v}_{k}=t_{k}^{\alpha}\widetilde{v}_{k}(t_{k}\cdot)$, then $\overline{u}_{k}\rightarrow u$ and $\overline{v}_{k}\rightarrow v$ weakly in $H_{r}^{1}(\R^{2})$, as $k\rightarrow\infty$. Next we show $u\in S_{1}$, $v\in S_{\omega}$ and $\overline{u}_{k}\rightarrow u$, $\overline{v}_{k}\rightarrow v$ in $H_{r}^{1}(\R^{2})$, as $k\rightarrow\infty$. This will be a contradiction.

Since $\widetilde{u}_{k}\in H_{r}^{1}(\R^{2})$, there exist $U_{k}\in C_{0}(\R^{2})$ and $V_{k}\in C_{0}(\R^{2})$ such that
\begin{align*}
\int_{\R^{2}}|\nabla \widetilde{u}_{k}-U_{k}|^{2}dx<\varepsilon,\
\ \ \int_{\R^{2}}|\widetilde{u}_{k}-V_{k}|^{2}dx<\varepsilon.
\end{align*}
Therefore
\begin{align*}
\|\nabla&(\overline{u}_{k}-\widetilde{u}_{k})\|_{2}^{2}\\
&=\int_{\R^{2}}|\nabla(t_{k}^{\al}\widetilde{u}_{k}(t_{k}x)-\widetilde{u}_{k}(x))|^{2}dx\\
&\leq2\int_{\R^{2}}|\nabla(t_{k}^{\al}\widetilde{u}_{k}(t_{k}x))-U_{k}(x)|^{2}dx+2\int_{\R^{2}}|\nabla\widetilde{u}_{k}(x))-U_{k}(x)|^{2}dx\\
&=2\int_{\R^{2}}|t_{k}^{\al+1}\nabla(\widetilde{u}_{k}(t_{k}x))-U_{k}(x)|^{2}dx+2\int_{\R^{2}}|\nabla\widetilde{u}_{k}(x))-U_{k}(x)|^{2}dx\\
&\leq4t_{k}^{2\al+2}\int_{\R^{2}}|U_{k}(t_{k}x)-U_{k}(x)|^{2}dx+2|t_{k}^{\al+1}-1|^{2}\int_{\R^{2}}U_{k}^{2}(x)dx+(4t_{k}^{2\al}+2)\varepsilon\\
&=12\varepsilon,
\end{align*}
and
\begin{align*}
\|&\overline{u}_{k}-\widetilde{u}_{k}\|_{2}^{2}\\
&=\int_{\R^{2}}|t_{k}^{\al}\widetilde{u}_{k}(t_{k}x)-\widetilde{u}_{k}(x)|^{2}dx\\
&\leq2\int_{\R^{2}}|t_{k}^{\al}\widetilde{u}_{k}(t_{k}x)-V_{k}(x)|^{2}dx+2\int_{\R^{2}}|\widetilde{u}_{k}(x)-V_{k}(x)|^{2}dx\\
&\leq4t_{k}^{2\al}\int_{\R^{2}}|V_{k}(t_{k}x)-V_{k}(x)|^{2}dx+2|t_{k}^{\al}-1|^{2}\int_{\R^{2}}V_{k}^{2}(x)dx+(4t_{k}^{2\al-2}+2)\varepsilon\\
&=12\varepsilon.
\end{align*}
It follows that $\|\overline{u}_{k}-\widetilde{u}_{k}\|_{H_{r}^{1}(\R^2)}\rightarrow0$ and $\|\overline{v}_{k}-\widetilde{v}_{k}\|_{H_{r}^{1}(\R^2)}\rightarrow0$ as $k\rightarrow\infty$. So
\begin{align*}
I(\widetilde{u}_{k},\widetilde{v}_{k})=I(\overline{u}_{k},\overline{v}_{k})+o_{k}(1)\geq \widetilde{E}_{1}+\widetilde{E}_{\omega}+o_{k}(1).
\end{align*}
Recalling that $\limsup_{k\rightarrow\infty}I(\widetilde{u}_{k},\widetilde{v}_{k})\leq \widetilde{E}_{1}+\widetilde{E}_{\omega}$, we obtain
\begin{align*}
\lim_{k\rightarrow\infty}J_{1}(\overline{u}_{k})=\widetilde{E}_{1},\ \ \ \lim_{k\rightarrow\infty}J_{\omega}\overline{v}_{k}=\widetilde{E}_{\omega}.
\end{align*}
Arguing as in the proof of Lemma \ref{Lemma 4.1}, we deduce that $u\not\equiv0$. Thanks to the lower semicontinuity of norms,
\begin{align*}
J_{1}(u)\leq\liminf_{k\rightarrow\infty}J_{1}(\overline{u}_{k})=B_{1}
\end{align*}
If $J_{1}(u)=E_{1}$, which yields that $\overline{u}_{k}\rightarrow u$ strongly in $H_{r}^{1}(\R^{2})$ and $u\in S_{1}$.
If not, we have
\begin{align*}
\|u\|_{H_{r}^{1}(\R^{2})}<\liminf_{k\rightarrow\infty}\|\overline{u}_{k}\|_{H_{r}^{1}(\R^{2})}.
\end{align*}
It follows that $u\not\in\mathcal{M}_{1}$. Then there exists a unique $t_{0}\in(0,1)$ such that $J(t_{0}^{\alpha}u(t_{0}\cdot))=0$. Thus, we have
\begin{align*}
J_{1}(t_{0}^{\alpha}u(t_{0}\cdot))&<\mathop{\lim}\limits_{n\rightarrow\infty}\Big(\frac{t_{0}^{2\alpha}}{2}\|\nabla\overline{u}_{k}\|_{2}^{2}+\frac{t_{0}^{2(\alpha-1)}}{2}\|\overline{u}_{k}\|_{2}^{2}+\frac{t_{0}^{6\alpha-4}}{2}B(\overline{u}_{k})-\frac{t_{0}^{2p\alpha-2}}{2p}\|\overline{u}_{k}\|_{2p}^{2p}\Big).
\end{align*}
Since $J_{1}(t^{\alpha}\overline{u}_{k}(t\cdot))$ has the maximum value at $t=1$ for all $k$. It follows that
\begin{align*}
J_{1}(t_{0}^{\alpha}u(t_{0}\cdot))&<\mathop{\lim}\limits_{n\rightarrow\infty}\Big(\frac{1}{2}\|\nabla\overline{u}_{k}\|_{2}^{2}+\frac{1}{2}\|\overline{u}_{k}\|_{2}^{2}+\frac{1}{2}B(\overline{u}_{k})-\frac{1}{2p}\|\overline{u}_{k}\|_{2p}^{2p}\Big)\\
&=\widetilde{E}_{1},
\end{align*}
which is a contradiction. That is $u\in\M_{1}$ and $J_{1}(u)=\widetilde{E}_{1}$, which yields $u\in S_{1}$, and $\overline{u}_{k}\rightarrow u$ strongly in $H_{r}^{1}(\R^{2})$ as $k\rightarrow\infty$. Finally, we can similarly prove $v\in S_{\omega}$ and $\overline{v}_{k}\rightarrow v$ strongly in $H_{r}^{1}(\R^{2})$ as $k\rightarrow\infty$. By Step 1, we know that $\widetilde{u}_{k}\rightarrow u$ and $\widetilde{v}_{k}\rightarrow v$ strongly in $H_{r}^{1}(\R^{2})$, as $k\rightarrow\infty$. This is a contradiction with the fact that $\widetilde{u}_{k}\in H_{r}^{1}(\R^{2})\setminus(S_{1})^{\delta_{0}}$ or $\widetilde{v}_{k}\in H_{r}^{1}(\R^{2})\setminus(S_{\omega})^{\delta_{0}}$. This completes the proof.
\qed

\

\noindent{\it Proof of $(i)$ of Theorem \ref{Theorem 1.1.}.} We divide the proof into several steps.

{\bf Step 1.}
$\M_{b}$ is nonempty. For each $(u,v)\in E\setminus\{(0,0)\}$, $J(t^{\alpha}u(t\cdot),t^{\alpha}v(t\cdot))$ is of the form $at^{2\alpha}+bt^{2(\alpha-1)}+ct^{6\alpha-4}-dt^{2p\alpha-2}$, which is positive for small $t$ and negative for large $t$. Thus, there exists $\widetilde{t}_{uv}>0$ such that $J(\widetilde{t}_{uv}^{\alpha}u(\widetilde{t}_{uv}\cdot),\widetilde{t}_{uv}^{\alpha}v(\widetilde{t}_{uv}\cdot))=0$. Thus, $\M_{b}$ is not empty.

{\bf Step 2.}
$\M_{b}$ is bounded away form zero, $i.e.$ $(0,0)\not\in\partial\M_{b}$. For each $(u,v)\in\M_{b}$,
\begin{align}
\label{1}
\begin{split}
F(u,v)&=\frac{p}{p\alpha-1}\Big(\alpha a(u,v)+(\alpha-1)b(u,v)+(3\alpha-2)c(u,v)\Big)\\
&\geq\frac{p}{p\alpha-1}(\alpha-1)\|(u,v)\|_{E}^{2}.
\end{split}
\end{align}
By the Sobolev embedding theorem there exists a constant $C>0$ such that for any $(u,v)\in\M_{b}$, $\|(u,v)\|_{E}^{2p}\geq C\|(u,v)\|_{E}^{2}$. Therefore, $\|(u,v)\|_{E}\geq\rho>0$ and the conclusion holds.

{\bf Step 3.}
$c_{b}>0$. For each $(u,v)\in\M_{b}$, combining (\ref{1})
\begin{align*}
I(u,v)&=\frac{1}{2}\|(u,v)\|_{E}^{2}+\frac{1}{2}\Big(B(u)+B(v)\Big)-\frac{1}{2p}F(u,v)\\
&=\frac{1}{2}\|(u,v)\|_{E}^{2}+\frac{1}{2}\Big(B(u)+B(v)\Big)-\frac{1}{2(p\alpha-1)}\Big(\alpha a(u,v)+(\alpha-1)b(u,v)+(3\alpha-2)c(u,v)\Big)\\
&=\Big(\frac{1}{2}-\frac{\alpha}{2(p\alpha-1)}\Big)a(u,v)+\Big(\frac{1}{2}-\frac{\alpha-1}{2(p\alpha-1)}\Big)b(u,v)+\Big(\frac{1}{2}-\frac{3\alpha-2}{2(p\alpha-1)}\Big)c(u,v)\\
&\geq\Big(\frac{1}{2}-\frac{\alpha}{2(p\alpha-1)}\Big)\|(u,v)\|_{E}^{2}.
\end{align*}
Then taking into account Step 2 and $p\alpha-1>\alpha$, we can obtain $c_{b}>0$.

{\bf Step 4:}
If $\{(u_{n},v_{n})\}$ is a minimizing sequence for $I$ on $\M_{b}$, then it is bounded. Let $\{(u_{n},v_{n})\}\in\M_{b}$ such that $I(u_{n},v_{n})\rightarrow c_{b}$.
As in Step 3, we get
\begin{align*}
I(u_{n},v_{n})=\Big(\frac{1}{2}-\frac{\alpha}{2(p\alpha-1)}\Big)a(u_{n},v_{n})+\Big(\frac{1}{2}-\frac{\alpha-1}{2(p\alpha-1)}\Big)b(u_{n},v_{n})+\Big(\frac{1}{2}-\frac{3\alpha-2}{2(p\alpha-1)}\Big)c(u_{n},v_{n}).
\end{align*}
Since the coefficients of $a(u_{n},v_{n})$, $b(u_{n},v_{n})$ and $c(u_{n},v_{n})$ are positive, then
\begin{align*}
I(u_{n},v_{n})\geq\Big(\frac{1}{2}-\frac{\alpha}{2(p\alpha-1)}\Big)\|(u_{n},v_{n})\|_{E}^{2},
\end{align*}
it follows that $\{(u_{n},v_{n})\}$ is bounded in $E$. Thus, there exists $(u,v)\in E$ such that, up to a subsequence
\begin{align*}
&(u_{n},v_{n})\rightarrow(u,v)\ \mbox{weakly\ in}\ E,\\
&(u_{n},v_{n})\rightarrow(u,v)\ \mbox{strongly\ in}\ L^{2p}(\R^{2})\times L^{2p}(\R^{2})\ \mbox{for}\ p\in(2,3],\\
&(u_{n},v_{n})\rightarrow(u,v)\ \mbox{a.e.}\ \mbox{in}\ \R^{2}.
\end{align*}
If $a(u,v)+b(u,v)=\lim\inf_{n\rightarrow\infty}a(u_{n},v_{n})+b(u_{n},v_{n})$, then it follows $(u_{n},v_{n})\rightarrow(u,v)$ strongly in $E$ as $n\rightarrow\infty$ and $(u,v)\neq(0,0)$, then $c_{b}$ is attained by $(u,v)$.\\
If $a(u,v)+b(u,v)<\lim\inf_{n\rightarrow\infty}a(u_{n},v_{n})+b(u_{n},v_{n})$, by the Lemma \ref{Lemma 2.1} and $J(u_{n},v_{n})=0$, we deduce that $J(u,v)<0$, then it follows that $(u,v)\not\in\M_{b}$ and $(u,v)\neq(0,0)$. Then there exists a unique $t_{0}\in(0,1)$ such that $J(t_{0}^{\alpha}u(t_{0}\cdot),t_{0}^{\alpha}v(t_{0}\cdot))=0$. Thus, we have
\begin{align*}
I(t_{0}^{\alpha}u(t_{0}\cdot),t_{0}^{\alpha}v(t_{0}\cdot))&<\mathop{\lim}\limits_{n\rightarrow\infty}\Big(\frac{t_{0}^{2\alpha}}{2}a(u_{n},v_{n})+\frac{t_{0}^{2(\alpha-1)}}{2}b(u_{n},v_{n})+\frac{t_{0}^{6\alpha-4}}{2}c(u_{n},v_{n})-\frac{t_{0}^{2p\alpha-2}}{2p}F(u_{n},v_{n})\Big).
\end{align*}
Since $J(t^{\alpha}u_{n}(t\cdot),t^{\alpha}v_{n}(t\cdot))$ has the maximum value at $t=1$ for all $n$.It follows that
\begin{align*}
I(t_{0}^{\alpha}u(t_{0}\cdot),t_{0}^{\alpha}v(t_{0}\cdot))&<\mathop{\lim}\limits_{n\rightarrow\infty}\Big(\frac{1}{2}a(u_{n},v_{n})+\frac{1}{2}b(u_{n},v_{n})+\frac{1}{2}c(u_{n},v_{n})-\frac{1}{2p}F(u_{n},v_{n})\Big)\\
&=c_{b},
\end{align*}
which is a contradiction.

{\bf Step 5.}
The minimizer $(u,v)$ is a regular point of $\M_{b}$, $i.e.$ $J'(u,v)\neq0$. To the contrary, suppose that $J'(u,v)=0$. For $(u_{t},v_{t})=(t^{\alpha}u(tx),t^{\alpha}v(tx))$, one has
\begin{align*}
&J(u_{t},v_{t})=t\frac{d}{dt}I(u_{t},v_{t}),\\
&\frac{d}{dt}J(u_{t},v_{t})=\frac{d}{dt}I(u_{t},v_{t})+t\cdot\frac{d^{2}}{dt^{2}}I(u_{t},v_{t}).
\end{align*}
Since $\frac{d}{dt}\Big|_{t=1}J(u_{t},v_{t})=0$, it follows that
\begin{align*}
2\alpha^{2}a(u,v)+2(\alpha-1)^{2}b(u,v)+2(3\alpha-2)^{2}c(u,v)-\frac{2(p\alpha-1)^{2}}{p}F(u,v)=0.
\end{align*}
Then, combining with $J(u,v)=0$, we get
\begin{align*}
0=&(\alpha^{2}-\alpha(p\alpha-1))a(u,v)+((\alpha-1)^{2}-(\alpha-1)(p\alpha-1))b(u,v)\\
&+((3\alpha-2)^{2}-(3\alpha-2)(p\alpha-1))c(u,v)
\end{align*}
The coefficients of $a(u,v)$, $b(u,v)$ and $c(u,v)$ in the above identity are negative, which is a contradiction.

{\bf Step 6.}
$I'(u,v)=0$. Thanks to Lagrange multiplier rule,there exists $\mu\in\R$ such that
\begin{align}
\label{2}
I'(u,v)=\mu J'(u,v).
\end{align}
We claim $\mu=0$.
There holds
\begin{align*}
\begin{cases}
\alpha a+(\alpha-1)b+(3\alpha-2)c-\frac{p\alpha-1}{p}d=0;\\
(1-2\alpha\mu)a+(1-2\mu(\alpha-1))b+3\big(1-\mu(6\alpha-4)\big)c-\big(1-\mu(2p\alpha-2)\big)d=0;\\
\big(2\mu(\alpha-1)-1\big)b+2\big(\mu(6\alpha-4)-1\big)c-\frac{\mu(2p\alpha-2)-1}{p}d=0.
\end{cases}
\end{align*}
The first equation holds since $J(u,v)=0$. The second one follows by multiplying (\ref{2}) by $(u,v)$ and integrating. The third one comes from Pohozaev equality.
We get
\begin{align*}
0=&\mu\Big(\big((2\alpha^{2}-(2p\alpha-2)\alpha\big)a(u,v)+\big(2(\alpha-1)^{2}-2(p\alpha-1)(\alpha-1)\big)b(u,v)\\
&+(2(3\alpha-2)^{2}-2(p\alpha-1)(3\alpha-2)c(u,v)\Big).
\end{align*}
All coefficients of $a(u,v)$, $b(u,v)$, $c(u,v)$ in the above identity are negative. This implies that $\mu=0$.

{\bf Step 7.} Thanks to Lemma \ref{Lemma 3.8}, the minimization $\{(u_{n},v_{n})\}$ can be chosen in $(S_{1})^{\delta}\times(S_{\omega})^{\delta}$, where $\delta>0$ is small such that $0\not\in S_{1}$ and $0\not\in S_{\omega}$. Hence, $u\not\equiv0$ and $v\not\equiv0$, $(u,v)$ is a nontrivial ground state solution of \eqref{question1}. In fact, since $(|u|,|v|)\in \N$ and $c_{b}=I(|u|,|v|)$, we conclude that $(|u|,|v|)$ is a nonnegative solution of (\ref{question1}). Using the strong maximum principle, we infer that $|u|,|v|>0$. Thus $(|u|,|v|)$ is a positive ground state solution of (\ref{question1}). This completes the proof.
\qed



\begin{thebibliography}{9999}
\bibitem{A-E} A. Ambrosetti, and E. Colorado, Bound and ground states of coupled nonlinear Schr\"odinger equations, {\it C.R. Math.,} {\bf 342}(2006), 453-458.

\bibitem{AE} A. Ambrosetti, and E. Colorado, Standing waves of some coupled nonlinear Schr\"odinger equations, {\it J. London Math. Soc.,} {\bf 75}(2007), 67-82.

\bibitem{AA} A. Azzollini and A. Pomponio,  Positive energy static solutions for the Chern-Simons-Schr\"odinger system under a large-distance fall-off requirement on the gauge potentials, {\it Calc. Var. Partial Differ. Equ.,} {\bf 60}(2021), 1-30.


\bibitem{B-W} T. Bartsch, and Z.-Q. Wang, Note on ground states of nonlinear Schr\"odinger systems, {\it J. Partial Differ. Equ.,} {\bf 19}(2006), 200-207.

\bibitem{J-H-J}	J. Byeon, H. Huh, and J. Seok, Standing waves of nonlinear Schr\"odinger equations with the gauge field, {\it J. Funct. Anal.,} {\bf 263}(2012), 1575-1608.

\bibitem{JHJ} J. Byeon, H. Huh, J. Seok, On standing waves with a vortex point of order $N$ for the nonlinear Chern-Simons-Schr\"odinger
equations, {\it J. Differ. Equat.,} {\bf 261}(2016) 1285-1316.

\bibitem{DT} D. Cassani, H. Tavares and J.-J. Zhang, Bose fluids and positive solutions to weakly coupled systems with critical growth in dimension two, {\it J. Differ. Equat.,} {\bf 269}(2020) 2328-2385.

\bibitem{G} G. V. Dunne, Self-Dual Chern-Simons Theories, Springer, 1995.

\bibitem{GT} D. Gilbarg and N. Trudinger, Elliptic Partial Differential Equations of Second Order, second edition, Grundlehren Math. Wiss., vol. 224, Springer, Berlin, 1983.

\bibitem{CR} C. R. Hagen, A new gauge theory without an elementary photon, {\it Ann. Phys.,} {\bf 157}(1984), 342-359.

\bibitem{C} C. R. Hagen, Rotational anomalies without anyons, {\it Phys. Rev.,} {\bf 31}(1985), 2135-2136.

\bibitem{H-H} H. Huh, Standing waves of the Schr\"odinger equation coupled with the Chern-Simons gauge field, {\it J. Math. Phys.,} {\bf 53}(2012).

\bibitem{R} R. Jackiw and S.-Y. Pi, Classical and quantal nonrelativistic Chern-Simons theory, {\it Phys. Rev. D.,} {\bf 42}(1990), 3500-3513.

\bibitem{S} R. Jackiw and S.-Y. Pi, Soliton solutions to the guaged nonlinear Schr\"odinger equations on the plane, {\it Phys. Rev. Lett.,} {\bf 64}(1990), 2969-2972.

\bibitem{C-F} C. Ji, and F. Fang, Standing waves for the Chern-Simons-Schr\"odinger equation with critical exponential growth, {\it J. Math. Anal. Appl.,} {\bf 450}(2017), 578-591.

\bibitem{Lin} T.-C. Lin, and J.-C. Wei, Ground state of N coupled nonlinear Schr\"odinger equations in $\mathbb{R}^{n}$, $n\leq3$,
{\it Comm. Math. Phys.,} {\bf 255}(2005), 629-653.

\bibitem{zzz} Z.-S. Liu, Z. Ouyang, and J.-J. Zhang, Existence and multiplicity of sign-changing standing waves for a gauged nonlinear Schr\"odinger equation in $\R^{2}$, {\it Nonlinearity,} {\bf 32}(2019), 3082-3111.

\bibitem{M-M-P} L. A. Maia, E. Montefusco, and B. Pellacci, Positive solutions for a weakly coupled nonlinear Schr\"odinger system, {\it J. Differ. Equat.,} {\bf 229}(2006).

\bibitem{P-R} A. Pomponio, and D. Ruiz, A variational analysis of a gauged nonlinear Schr\"odinger equation, {\it J. Eur. Math. Soc.,} {\bf 17}(2015), 1463-1486.

\bibitem{AD} A. Pomponio, and D. Ruiz, Boundary concentration of a gauged nonlinear Schr\"odinger equation on large balls, {\it Calc. Var.
Partial Differ. Equ.,} {\bf 53}(2015) 289-316.

\bibitem{Rab} P. H. Rabinowitz, On a class of nonlinear Schr\"odinger equations, {\it Z. Angew. Math. Phys.,} {\bf 43}, 270-291(1992).

\bibitem{RD} D. Ruiz, The Schr\"odinger-Poisson equation under the effect of a nonlinear local term, {\it J. Funct.Anal.,} {\bf 237}(2006), 655-674.

\bibitem{Sirakov}  B. Sirakov, Least energy solitary waves for a system of nonlinear Schr\"odinger equations in $\mathbb{R}^{n}$, {\it Comm. Math. Phys.,} {\bf 271}(2007), 199-221.

\bibitem{Willem} M. Willem, Minimax Theorems, Birkh\"auser, Boston, 1996.

\bibitem{Jian} J. Zhang, W. Zhang, and X. Xie, Infinitely many solutions for a gauged nonlinear Schr\"odinger equation, {\it Appl. Math. Lett.,} {\bf 88}(2019), 21-27.


\end{thebibliography}
\end{document}